\documentclass[12pt,draft]{amsart}
\usepackage{amssymb,a4,eucal}
                                                  
\def\C{\mathcal{C}}

\def\I{\mathcal{I}}
\def\E{\mathcal{E}}
\def\cqfd{\qed}

\def\un{\mbox{$1\!\!1$}}

\def\dint{\displaystyle \int }

\DeclareSymbolFont{AMSa}{U}{msa}{m}{n}
\DeclareSymbolFont{AMSb}{U}{msb}{m}{n}
\DeclareMathSizes{12}{11}{7}{5}
\DeclareMathSizes{10}{10}{7}{5}
\DeclareMathSymbol{\Rset}         {\mathalpha}{AMSb}{"52}
\DeclareMathSymbol{\Zset}         {\mathalpha}{AMSb}{"5A}
\DeclareMathSymbol{\Nset}         {\mathalpha}{AMSb}{"4E}
\DeclareMathSymbol{\Qset}         {\mathalpha}{AMSb}{"51}
\DeclareMathSymbol{\Cset}         {\mathalpha}{AMSb}{"43}
\DeclareMathSymbol{\Bset}         {\mathalpha}{AMSb}{"42}
\DeclareMathSymbol{\Kset}         {\mathalpha}{AMSb}{"4B}

\newtheorem{thm}{Theorem}
\newtheorem{cor}[thm]{Corollary}
\newtheorem{lem}[thm]{Lemma}
\newtheorem{prop}[thm]{Proposition}

\theoremstyle{definition}
\newtheorem{defn}{Definition}
\newtheorem{rk}{Remark}

\theoremstyle{remark}
\newtheorem{ex}{Example}

\def\supp{\mathop{\textup{supp\,}}}

\def\dist{\mathop{\textup{dist\,}}}

\def\dint{\displaystyle \int }

 \author{Youssef Jabri \and  Mimoun Moussaoui}
 \address{Department of mathematics,
         University Mohamed I,  Box 524,
         60000 Oujda, Morocco }

\begin{document}
 \title{On the Linking Principle}
 \date{}
 
 \maketitle
   
 \section{Introduction}
 
 During the last twenty years, many minimax theorems that have proved 
 to be very useful tools in finding critical points of functionals have 
 been established.  
 They have all in common a geometric intersection property   
 known as the linking principle.
 
 Our purpose in this paper is to give a linking theorem that strengthens and 
 unifies some of these works. We think essentially to Ambrosetti-Rabinowitz 
``mountain pass theorem''~\cite{A-R}, Rabinowitz  
 ``multidimensional mountain pass theorem''~\cite{R1}, Rabinowitz ``saddle point 
 theorem''~\cite{R2}  and  Silva's variants of these results~\cite{silva}.
 
 We focus our attention especially on  ``the limiting case'', known 
to be true for the mountain pass principle \cite{R4}, where some information on 
 the location of the critical points is given. We give two forms of this theorem, in
the first part of the paper, the first one is established via a deformation lemma
and in the second part we use Ekeland's variational principle to get the second one.
Unfortunately, we could not remove  the finite dimension condition that 
appears in the results~\cite{R4} and~\cite{R3}. This finite dimension 
assumption is dropped only when dealing with a special kind of functionals 
having a predefined shape~\cite{B-R,Y,silva}.

\section{Using the deformation lemma}

 We start by giving a short description of the notion of linking that seems to be
 more natural than the initial one introduced by Benci and Rabinowitz. 
 It formalizes a topological property of intersection that appears in 
 all the results cited above.
 \begin{defn}
 Let S be a closed subset of a Banach space V, Q is a subset of V with 
 relative boundary $ \partial Q$, we say S and $ \partial Q$ link if : \\
 \textup{(L1)} $ \quad S \cap \partial Q = \varnothing$ \\
 \textup{(L2)} $ \quad$ for all $\gamma \in {\C}(V,V)$ such that 
 $ \gamma \bigm\vert_{\partial Q} = I_d$ 
 we have $ \gamma (Q)\cap S \neq \varnothing $ \\
 More generally, if $ \Gamma$ is a subset of $ {\C}(V,V)$, then $S$ and
 $ \partial Q$ are said to link with respect to $ \Gamma$ if the relation 
 \textup{(L1)} holds and \textup{(L2)}  is satisfied for any $ \gamma \in \Gamma$. 
 \end{defn}
 
 The following examples yield respectively the geometry of the saddle point 
 theorem and the generalized mountain pass theorem. 
 
 \begin{ex}
 	Let $V=V_1 \oplus V_2$ be a space decomposed into two closed subspaces 
 	 $V_1$ and $V_2$ with $\dim V_2 < \infty$.\\
 	 Let $S=V_1$ and $Q=B_R(0) \cap V_2$ with relative boundary
 	 $$ 
	 \partial Q =\big\{u \in V_2;\  ||u||  =R \big\}.
	 $$ 
 	 Then $S$ and $\partial Q$ link.
 \end{ex}
 \begin{ex}
	 Let $V=V_1 \oplus V_2$ as in Example~1 and let $e \in V_1$ with $||u||$ be given. 
	 Suppose $ 0 < \rho < R_1, 0 < R_2 $ and let
	 $$ 
	 S = \big\{u \in V_1;\ ||u|| = \rho \big\},
	 $$
	 $$
	 Q=\big\{se+u_2;\ 0\leq s \leq R_1,  u_2 \in V_2 \text{ and }  ||u_2|| \leq R_2 \big\}
	 $$
	 with relative boundary
	 $$ 
	 \partial Q=\big\{ se+u_2;\ s \in \{0,R_1\}\text{ or }||u_2||=R_2 \big\}.
	 $$
	 Then $S$ and $ \partial Q$ link.
\end{ex}

 These examples are now well known in the specialized literature (see 
for instance~\cite{B-R}). For a proof see \cite{St} or \cite{B-B-F}. 
A short hint is given in the Appendix.\\
 
In minimax theorems, the functional must verify some compactness property known
 as the Palais-Smale condition. We recall it. 
 \begin{defn}
 Let $\Phi \in {\C}^1(V,\Rset)$ and $c \in \Rset$. We say that $\Phi$ satisfies
 the  Palais-Smale Condition \textup{(P.S.)} if the existence of a 
 sequence $(u_n)_n$ in $V$ such that $\Phi (u_n)$ is bounded and $\Phi' (u_n) 
 \to 0$ as $n \to \infty$ implies that $(u_n)_n$ possesses a convergent subsequence.
 And $\Phi$ satisfies the local Palais-Smale condition 
 at $c$ denoted by  \textup{(P.S.)}$_c$, if the existence of a sequence $(u_n)_n$
 in $V$ such that $\Phi (u_n) \to c$ and $\Phi' (u_n) \to 0$ as $n \to \infty$ 
 implies that it is precompact.
 \end{defn}

 \subsection{A Generalized linking theorem}

 We will use in the sequel the following notations:
 $$
 \Phi^a=\big\{ v \in V ;\ \Phi (v) \leq a \big\}, 
 \qquad
 \Phi^a_b=\big\{ v \in V;\ b \leq \Phi (v) \leq a \big\},
 $$
 $$
 \mathcal{K}_c=\big\{ v \in V;\ \Phi' (v)=0 \text{ and } 
 \Phi(v)=c\big\},
 $$
 $$
N_\delta (E)=\big\{ v \in V;\ \dist(v,E)= ||v-E|| \leq\delta \big\},
 $$
 and
$$
\tilde V=\big\{  v\in V;\ \Phi'(v) \neq 0 \big\}.
$$
 
 Let us now state the abstract critical point theorem we announced in the 
 beginning.
 
 \begin{thm}\label{premier}
 Suppose that $\Phi \in {\C}^1(V,\Rset)$, $S \subset V$ is a closed subset and $Q 
 \subset V$ satisfy:\\
 \textup{(a)} $ S $ and $ \partial Q$ link with respect to $ \Gamma$. \\
 \textup{(b)} There exists $ \alpha \in \Rset$ such that
 $$
 \Phi \bigm\vert_{ \partial Q } \leq \alpha \leq \Phi \bigm\vert_S 
 $$ 
 \textup{(c)} There exists $ \gamma_0 \in { \C}(V,V)$ such that 
 $$
 \sup_{u \in Q} \Phi \left( \gamma_0(u) \right) \  < \infty. 
 $$  
 Let  
 $$
 c_\Gamma=\inf_{ \gamma \in \Gamma} \sup_{u \in Q} \Phi \left( \gamma (u)
\right),
$$ 
if $\Phi$ satisfies \textup{(P.S.)}$_{c_ \Gamma }$ then ${c_ \Gamma}$  
defines a critical value and $c_{\Gamma} \geq \alpha$. \\  
 Moreover if $c_{\Gamma} = \alpha$ then 
 $\mathcal{K}_{c_{\Gamma}} \cap S \neq \varnothing$. \\
 By $ \Gamma$ we denote one of the sets
 \begin{itemize} 
 \item $\Gamma_1 =\big\{\gamma \in {\C} (V,V);\ \gamma|_{ \partial Q}
 = I_d \big\}$ 
 \item $\Gamma_2 =\big\{ \gamma \in \mathcal{K} (V,V)  \colon \text{ 
 set of compact maps};\ \gamma|_{ \partial Q} = I_d \big\}$ 
 \item $ \Gamma_3 =\big\{ \gamma \in \mathcal{H} (V,V) \colon \textrm{set 
 of homeomorphisms};\  \gamma|_{ \partial Q} = I_d \big\}$ 
 \item $\Gamma_4 = \big\{ \gamma \in \mathcal{K}(V,V) \colon  \gamma 
 \text{(closed set)  =  closed  set  and  } \gamma I_{\partial Q} = I_d \big\}$
 \item $\Gamma_5 = \big\{ \gamma \in \Gamma_1\text{ such that } \gamma \text{(closed 
  bounded  set)   =  compact  set } \big\}$\\
 \end{itemize} 
 Compact means continuous and maps bounded sets into relatively compact 
 ones.\\
 \end{thm}

 \begin{rk}
 	\begin{itemize}
 	 \item The assumption (c) is satisfied when $Q$ is compact.
 	 \item In the classical results, $Q$ is compact and in general it is assumed that
 	 there  exists $ \beta < \alpha $ such that 
 	 $$
 	 \Phi |_{ \partial Q } 
 	 \leq \beta < \alpha \leq \Phi |_S. 
 	 $$
 	 \item If we denote by $\mathcal{F}$ the set of the $\Gamma_i$ defined 
 	 above, it is still possible to consider for $\Gamma$ any finite intersection 
 	 in $\mathcal{F}$. 
 	 \end{itemize} 
 \end{rk}
  
 \begin{rk}
	The set used in general in the classical results is $\Gamma_1$. The 
	set   
	$\Gamma_3$ has also been used, but $ \Gamma_2$, $\Gamma_4$  and $\Gamma_5$ 
	have never been used  at our best knowledge. They have a
	particular importance because  the sets in examples 1 and 2 link with 
	respect to $\Gamma_2$, $\Gamma_4$ and $\Gamma_5$ too,  \emph{even when 
	$\dim V_2 $ is infinite}\footnote{This could lead to removing the finite 
	dimension assumption of $V_2$}.  
	Unfortunately, if these sets are small enough to make that $\partial Q$ 
	and $S$ link, they are empty when $\dim V_2$ is infinite (See Remark~10 at the end).
\end{rk}

 Notice also that we do not require that the values of $\Phi$ on $\partial Q$ and $S$ are
\emph{strictly separated.}\\ 
 The reader is supposed familiarized with these classical results. So even 
 improved, we will not state them here.\\
 
 We want to point out also that\\ 
 $ \left\{~
 \begin{tabular}{l}
 $c_{\Gamma_3} \geq c_{\Gamma_1} \geq \alpha $, and \\
 $c_{\Gamma_5} \geq c_{\Gamma_4} \geq c_{\Gamma_2} \geq  c_{\Gamma_1} \geq \alpha $\cr  
 \end{tabular} 
 \right.
 \qquad$
where the $c_{\Gamma_i}$'s are defined in Theorem~1\\ 
 
 If $\Gamma_i\subset\Gamma_j$ then $c_{\Gamma_j} \geq c_{\Gamma_i} 
 \geq c_{\Gamma_1} \geq \alpha$ and there is no reason they must be equal 
 in general. Nevertheless this  has not been verified.
 
 The proof follows a classical method. We use a slight modification of 
 the standard deformation lemma. It can be found in~\cite{Y} when 
 $V = V_1 \oplus V_2$ is a Hilbert space and $\Phi(u)=(Lu,u)+b(x)$ with 
 $L = L_{1} P_{1} + L_{2} P_{2}$ where  
 $P_{i} \colon V \to V_{i}$ is the projection onto~$V_{i}$, the 
 operator $L_{i} \colon V_{i} \to V_{i}$ is bounded selfadjoint 
 and the derivative $b'$ of $b$ is compact.
 
 \begin{lem}\label{lemme1}
 Let $V$ a Banach space, $D$ and $E$ two closed subsets of $V$ and $f \colon
 V \to \Rset$ a ${\C}^1$-functional satisfying \textup{(P.S.)}$_c$ for 
 $c \in \Rset$ such that \\  
\textup{(L1)} $D \cap E = \varnothing \quad$  and $\quad E \cap 
 \mathcal{K}_c = \varnothing,  \quad$  and  \\ 
\textup{(L2)} $f |_E \leq c \leq f |_D $. \\
 Then, there exist $ \varepsilon > 0$ and $ \eta \in { \C}( \Rset \times V,V)$
 such that \\
 \begin{tabular}{rl}
 \textup{(i)}&$\eta_{t}= \eta (t,.) \colon V \to V$ is a homeomorphism
and $\eta_{t}^{-1} = \eta_{-t},$ for all t. \\
\textup{(ii)}&$ f( \eta(t,x)) \leq f(x)$, for all $t \geq 0$, \\
\textup{(iii)}&$\eta (t,x) = x$,  for all $t \in \Rset$, for any $x \in 
D$, and \\ 
\textup{(iv)}&$f( \eta (1,x)) \leq c- \varepsilon$, for all $x \in E$ \\
\end{tabular}
 \end{lem}
  There are many deformation lemmas in the literature,  some  were
 proved sometimes with weaker conditions than \textup{(P.S.)}
 (see for example~\cite{silva,B-B-F}). The proof of this one is 
 practically given in \cite{Y}, we will sketch it in the Appendix for the 
 convenience of the reader.
 
 \begin{proof}{Proof of Theorem~\ref{premier}}
 $\bullet$ Notice that by assumption (c), $c_{\Gamma} < \infty $. \\
 $\bullet$ since $S$ and $\partial Q$ link, clearly $c_{\Gamma} \geq \alpha$. 
 So $c_{\Gamma}$ is well defined.
 
 We will distinguish the two cases:
 \begin{itemize}
 \item \textbf{Suppose that} $\mathbf{c_{\Gamma} > \alpha }$.
 Let us suppose that $\mathcal{K}_{c_{\Gamma}} = \varnothing$. Set 
 $\overline{\varepsilon} =(c_{\Gamma}-\alpha)/2$. By the classical deformation lemma 
 (see \cite{R3}) we have that for any $\varepsilon < \overline{\varepsilon}$, 
 there  exists $\eta\in{\C}\left([0,1] \times V,V\right)$ such that
$$
 \eta \left( 1,\Phi^{c_{\Gamma} + \varepsilon} \right) \subset 
 \Phi^{c_{\Gamma} - \varepsilon}. 
$$
 But by the definition of $c_{\Gamma}$, there exists $\gamma \in \Gamma$ such
that $ \sup _ {u \in Q} \Phi(\gamma(u))<{c_{\Gamma} + \varepsilon} $. 
Setting $\gamma' = \eta(1,.)\,\circ\, \gamma$ we have  $\sup _ {u \in Q} \Phi
(\gamma'(u)) < {c_{\Gamma} - \varepsilon}$. This implies that $\gamma' |_{\partial
Q} = I_d$, so $\gamma' \in \Gamma$ but this contradicts the definition of
$c_{\Gamma}$. 
 \item Suppose $\mathbf{c_{\Gamma}=\alpha}$.
We claim that $\mathcal{K}_{c_{ \Gamma}} \cap S \neq \varnothing$. 
Indeed, if by contradiction this was not the case. Then, since
$$
 -\Phi |_S \leq -c_{\Gamma} \leq -\Phi |_{\partial Q} 
$$
 By Lemma 2.1, there exists $\varepsilon$ and $\eta \in {\C}(\Rset 
 \times V,V)$ such that $\eta (t,x) = x$ on $\partial Q$ and $-\Phi \left(\eta 
 (1,x) \right) \leq -c_{\Gamma}-\varepsilon$ for all $x$ in $S$. But the definition
of $c_{\Gamma}$ implies  that there exists $\gamma \in \Gamma$ such that  
$$
 \Phi(\gamma(x)) < c_{\Gamma}+\varepsilon \quad text{for all } x \in Q.\eqno(*) 
$$ 
 Let $\gamma'(x)=\eta_{1}^{-1}(\gamma(x))$, then by (iii) and the choice 
 of $\Gamma_i$, the functional $\gamma' \in \Gamma$. \\ 
 Now, since $S$ and $\partial Q$ link, there exists $\bar x \in Q$ such that 
 $\gamma'(\overline{x}) \in S$ and then $\Phi(\gamma(\overline{x}))= 
 \Phi \left(\eta(1,\gamma'(\overline{x}))\right) 
 \geq c_{\Gamma}+\varepsilon$ and this contradicts ($*$). 
 \end{itemize}
 \end{proof}
 
 In \cite{silva}, Silva has shown  with an other deformation lemma some variants
 of the classical results cited above, which, at least at a first look, don't 
 seem to be a part of the general linking principle stated here.  
 We mention a saddle point theorem with some \emph{global estimates}.
 \begin{thm} 
  \label{silva1}
  Let $V=V_1 \oplus V_2$ be a real Banach space with $\dim V_2<+\infty$,  $\Phi\in
 {\C}^1(V,\Rset)$ such that
 $$
 \Phi |_{V_2} \leq \beta, \qquad \beta \in \Rset. \leqno(a)
 $$
 And
 $$ 
 \Phi |_{V_1} \geq \alpha, \qquad \alpha \in \Rset. \leqno(b)
 $$
 If $\Phi$ verifies \textup{(P.S.)}$_c $, for any $c \in [\alpha,\beta]$, then 
 $\Phi$ admits a critical value $c_{\Gamma} \in [\alpha,\beta]$, characterized by a 
 minimax argument.
 \end{thm}
 This result can also be obtained with Theorem 2.1 even if it uses  global 
 estimates instead of local ones in which the linking sets are explicited. 
 Indeed, by Silva's deformation lemma, there exist $\varepsilon,\ R_0 >0$ and
 $\eta \in {\C}([0,1] \times V,V)$ \   such that 
 $$
 \eta (1,u) \in \Phi^{\alpha - \varepsilon }, \quad \text{for all } u \in \Phi^ 
 \beta \setminus \Bset(0,{R_0})
 $$ 
 Since by Example 1, $S=V_1$ and $\partial Q $ link when $Q=\Bset(0,{R_0})\cap
V_2$. We will be in the situation of Theorem 1 if we prove that $c=c_{\Gamma}
\in [\alpha,\beta]$. And this last statement is  true. Indeed, since $\partial Q$ and
 $S$ link. For all $\gamma \in  \Gamma$, we have $\gamma (Q)\cap S \neq \varnothing,$ so
$\sup_Q \Phi \circ \gamma \geq \inf_S \Phi \geq \alpha$ and hence 
$$
c_{\Gamma} = \inf_{\gamma \in \Gamma} \sup_{ Q } \Phi
\circ \gamma \geq \alpha.
$$
The fact that the restriction $\gamma |_{\partial Q}=I_{d_{V_2}}$ implies that
$P_2 \circ  \gamma  |_{\partial Q} = I_{d_{V_2}}$ where $P_2 \colon V \to
V_2$ is the projection onto $V_2$. But  $\gamma' = P_2 \circ \gamma \in
\Gamma$ and $\gamma' (Q) = P_2 \circ  \gamma (Q) \subset V_2$. Therefore
$\sup_Q \Phi \circ \gamma' \leq \beta$ and  $c_\Gamma = \inf_{\gamma \in
\Gamma} \sup_{Q} \Phi \circ \gamma~\in  [\alpha,\beta]$ .\\
 
 We mention equally an other abstract critical point theorem by Silva:
 
 \begin{thm}[Silva]
 \label{silva2}
 Let $V = V_1 \oplus \Rset e \oplus V_2$ be a real Banach space with 
$\dim V_2<\infty$, where $|| e || =1$ and $\Phi \in {\C}^1(V,\Rset)$ 
are such that: 
$$
\text{ There exists }\rho >0\text{ and } \alpha \geq 0\text{ such  that } 
  \Phi |_{S_{\rho}} \geq \alpha.
   \leqno(a)
$$
 where $ S_{\rho} = {\overline{\left( V_1 \setminus \Bset(0,\rho) \right)}
 \cup \left( ( \Rset^+e \oplus V_1) \cap \partial \Bset(0,\rho) \right)}$.
 $$
 \Phi |_{V_2 \oplus \Rset e} \leq \beta. \leqno(b)
 $$
 If $\Phi$ satisfies \textup{(P.S.)}$_c$ for all $c$ in $\left[\alpha,\beta \right],$ 
 then $\Phi$ possesses a critical value $c_\Gamma$ in $\left[ \alpha,\beta \right]$  
 characterized by a minimax argument.
 \end{thm}

 \begin{rk}
	 We point out that in Theorem~\ref{deuxieme} we have $\alpha \leq \Phi(0) \leq \beta$ 
	 and $\alpha \leq \Phi(\rho e) \leq \beta $  in Theorem~\ref{silva2}.
 
 In this last result, Silva's deformation lemma implies that (see \cite{silva}):\\ 
 There exists $R_0 > \rho$ a real $ \varepsilon \in \left] 0,\alpha 
 \right[$ and a deformation $\eta
 \in {\C} \left( [0,1] \times V,V \right)$ such that $\ \eta (1,u) \in \Phi^{\alpha -
 \varepsilon }$ for all $u \in \left(\Phi^\beta\setminus \Bset(0,{R_0})\right)\cap(V_2
 \oplus \Rset e)$. But if we set {$Q=\Bset(0,{R_0})\cap( V_2 \oplus \Rset e )$} then $S_{\rho}$,
 as defined  in Theorem~3, and $\partial Q$  link. We give the proof in 
 the Appendix.
 There too, Theorem~\ref{premier} applies if we prove that $c_\Gamma \in [\alpha,\beta]$. 
 And this is proved similarly to what we did in the proof of Theorem~\ref{silva1}. 
\end{rk}
 \begin{rk}
	 In all the results cited here, no finite dimension is required to
	prove the linking with respect to the particular sets $\Gamma_i$ in use.  
\end{rk}

\begin{rk}
	The theorems seen here are of course applied in finding stationary points 
	 of functionals and then to solve variational problems. But  
	 unfortunately they have the drawback to impose  
	 the compactness condition {(P.S.)}$_c$ which is  unlikely to verify 
	 and seems rather restrictive. But recently Struwe mentioned in \cite{St}
	 that the failure of {(P.S.)}$_c$ at certain levels reflects, using
	 physicists terminology, phenomena related to ``phase transition'' or  ``particle
	creation'' at these levels, (for more details see \cite{St}).  
\end{rk}

\section{Using Ekeland's variational principle}

 Now, we will use some convex analysis results and Ekeland's
variational principle to prove an other variant of the abstract theorem stated in the
first part. The fact that the critical point is located on $S$ in the limiting case, 
will be confirmed again.

\subsection{Ekeland's variational principle}

When a functional is l{.}s{.}c{.} (lower semi-continuous) in a reflexive Banach space, it
possesses a minimum if and only if it has a bounded minimizing sequence. But if
the space is not reflexive or the functional is only lower semi-continuous, 
we can say no thing similar. Ekeland has introduced in 1972 a ``variational 
principle'' that applies well in such situations.
\begin{thm}[Ekeland's variational principle] Let $(X,d)$ be a complete metric space
and $\Phi\colon X\to \Rset\cup \{+\infty\}$ a lower semi-continuous functional,
bounded from below and not identical to $+\infty$. Then, for all $\varepsilon>0$,
each $\delta>0$ and each $x\in X$ such that
$$
\Phi(x)\leq \inf_{x\in X}\Phi(x) +\varepsilon.
$$
There exists $y\in X$ with the properties\\
a) $\Phi(y)\leq  \Phi(x)$\\
b) $\dist(x,y)\leq  \delta$\\
c) $\Phi(z)>\Phi(y)-\dfrac{\varepsilon}{\delta}\,\dist(z,y)$ for all $z\neq y$ in $X$
\end{thm}
This principle has been successfully used many times to prove the existence 
of ``almost critical points'' of unbounded functionals via a minimax method.
\begin{thm}
\label{strict}
Let $X$ be a Banach space, $K$ a compact metric space, $K_0\subset K$ a closed
subset, $\mu \in \C(K_0,X)$ and $E=E(K,K_0,X,\mu )$ the complete metric space defined
by
$$
E=\big\{ m\in \C(K,X);\ m|_{K_0}=\mu \big\}
$$
endowed with the distance of uniform convergence on $K$.

Let $f\in \C^1(X,\Rset)$ and set
$$
c=\inf_{m\in E}\max_{s\in K} f(m(s))
$$
and 
$$
d=\max_{s\in K_0}f(\mu (s)).
$$
Suppose that 
$$
c>d.
$$
Then, for each $\varepsilon \in ]0,c-d[$ and each $p\in E$ such that
$$
\max_{s\in K}f(p(s))\leq  c+\varepsilon
$$
there exists $u\in X$ such that
$$
c-\varepsilon \leq  f(u)\leq  \max_{s\in K} f(p(s)),
$$
$$
\dist(u,p(K))\leq  \sqrt {\varepsilon}
$$
and
$$
\|f'(u)\|\leq  \sqrt {\varepsilon}.
$$
\end{thm}
In \cite{Ma}, \cite{Defi} and \cite{M-W}, two nice proofs of this result using Ekeland's
variational principle are given. 

This theorem contains as particular cases the theorems cited in the first part.
Unfortunately, in this theorem the notion of ``linking'' described before is hidden,
and the strict separation of the values of the functional on $\partial Q$ and $S$ is
essential in the proofs.\footnote{In fact, the linking is reintroduced 
again by the authors to obtain the particular cases cited above.}

We will use Ekeland's principle to prove an abstract critical point theorem 
close to Theorem~\ref{premier} 
where we both exhibit the notion of linking and treat the ``limiting case''. We want
to point out that, in this part, in addition to~\cite{Ma,Defi,M-W} we have
have been inspired by the paper of~\cite{G-P} devoted to the particular case of the
mountain pass theorem.
\begin{thm}
\label{deuxieme}
Let $X$ be a Banach space, $Q$ a compact subset and $S$ a closed subset of $X$
such that $S$ and $\partial Q$ link. 

Let $f\colon X\to \Rset$ be a $\C^1$-functional and set
$$
E=\{ \gamma\in \C(Q,X) ;\ \gamma|_{\partial Q} =I_d \},
$$
$$
c=\inf_{\gamma\in E}\max_{x\in Q} f(\gamma(x))
$$
and
$$
c_0=\inf_S f\text{, and }\qquad d=\max_{\partial Q}f
$$
Suppose that 
$$
c_0\geq  d.\eqno(*)
$$
Then $f$ admits a sequence of {``almost critical''} points $(u_n)_n\subset E$ such
that $f(u_n)\to c$ and $f'(u_n)\to 0$. In the  case $c=c_0$, (in
this case $c=c_0=d$) we have some informations on the location of these points.
For each $0<\varepsilon<\max\{1,\dist(\partial Q,S)\}/{2}$, there exists
$x_\varepsilon\in X$ such that:\\
{\rm (i)} $c\leq  f(x_\varepsilon \leq  c+\dfrac{5}{4}\varepsilon^2$,\\
{\rm (ii)} $\dist(x_\varepsilon, S)\leq  \dfrac{3}{2}\varepsilon$, and \\
{\rm (iii)} $\|f'(x_\varepsilon)\|\leq  \dfrac{\varepsilon}{2}.$
\end{thm}
\textsc{Proof.} Since $S$ and $\partial Q$ link, the relation ($*$) implies that
$$
c\geq  c_0\geq  d.
$$
If $c>c_0$ or $c_0>d$ then $c>d$ and we are in the situation of
Theorem~\ref{strict}. So, it suffices to treat the limiting case $c=c_0=d$. In fact,
instead of ($*$) we will use only the condition $c=c_0$, but in this case 
it is legitimate to wonder if we have not the situation
$$
d>c_0.
$$
The answer is negative, we have really $c_0\geq  d$ and therefore ($*$). Indeed, if
there exists $t\in \partial Q$ with $f(t)>c_0$. Since $t\in \partial Q$, we would have
$\gamma(t)=t$ for each $\gamma \in E$ and hence $f(\gamma(t))=f(t)>c_0$ for
each $\gamma\in E$.

Therefore
$$
c\geq  f(t)\geq  c_0,
$$
which contradicts the equality $c=c_0$.

Before continuing the proof, we will recall some known results uses in 
the sequel and whose proofs can be found in standard convex analysis literature.

\subsection*{Preliminaries}

Let $\Phi\colon X\to \Rset$ be a functional, the subdifferential of $\Phi$ is the
multifunction $\partial \Phi\colon X\to \mathcal{P}(X^*)=2^{X^*}$ defined by
$$
\partial \Phi(x)=\big\{ \mu \in X^*;\ \Phi(y)\geq  \Phi(x)+\langle \mu, 
y-x\rangle, ~\forall y\in X\big\}  
$$
where $X^*$ is the dual space of $X$.
\begin{prop}
Let  $\Phi\colon X\to \Rset$ be continuous and convex. Then for each $x,y\in X$
$$
\lim_{t\downarrow0}\dfrac{\Phi(x+ty)-\Phi(x)}{t}=\max_{\mu \in \partial 
Q(x)}\langle \mu ,y  \rangle.
$$
\end{prop}
\begin{prop}
Let  $\Phi\colon X\to \Rset$ be continuous and convex. Then $\partial \Phi(x)$ is non-empty,
convex and $w^*$-compact in $X^*$ for each $x\in X$.
\end{prop}

Let $K$ be a compact metric space and $\E=\C(K,\Rset)$ the Banach space of
continuous functions for the norm $\|x\|=\max\{x(t);\ t\in K\}$. By Riesz
representation Theorem, the dual space $\E^*$ of $\E$ is isometric isomorphic to the
Banach space of regular Radon measures defined on the $\sigma$-algebra of Borel
sets of $K$.

Recall that:\\
$\bullet$ A Radon measure $\mu $ is positive ($\mu \geq  0$) if 
$\langle\mu ,x\rangle \geq  0$ for all $x\in E$ such that $x(t)\geq  0$ 
for each $t\in K$.\\  
$\bullet$ A Radon measure has mass one if $\langle \mu ,\un\rangle=1$, where $\un\colon
K \to\Rset$ is the constant function $\un(t)=1$ for any $t\in K$.\\
$\bullet$ A Radon measure vanishes in an open subset $U\subset K$ if $\langle
\mu ,x\rangle =0$ for each $x\in \E$ such that the support of $x$ is a compact set
$K\subset U$. If $\mu $ vanishes in a collection of open subsets $(U_\alpha )_\alpha$,
then $\mu $ vanishes in $\cup_\alpha U_\alpha$, therefore there exists a largest
open set $\tilde U$ where it vanishes.
The support of $\mu $, denoted by $\supp\mu $, is defined by $\supp\mu =K\setminus \tilde
U$.
\begin{prop}
Let $\theta\colon\to \Rset$, $x\mapsto \theta (x)=\max\{x(t);\ t\in K\}$. Then 
$\theta$ is convex and continuous. Moreover, for each $x\in E$
$$
\mu\in \partial \theta(x)\iff\left\{
\begin{array}{l}
\mu\geq 0, \langle\mu,\un\rangle=1,\\
\supp\mu \subset\big\{ t\in K;\ x(t)=\theta(x)\big\}
\end{array}
\right.
$$
\end{prop}
\textsc{The proof of theorem \ref{deuxieme} continued:}\\
We recall that we are in the limiting case $c=c_0$.\\ 
Let $g\in E$ such that
\begin{equation}
\max_{t\in Q}f(g(t))<c+\dfrac{\varepsilon^2}{4}
\end{equation}
Let us denote by
$$
S_{g,\varepsilon}=\big\{ x\in X;\ \dist(g(x),S)<\varepsilon\big\}
$$
and set
$$
A=Q\cap S_{g,\varepsilon} =\big\{x\in Q;\ \dist(g(x),S)<\varepsilon\big\}
$$
$A$ is non-empty because by the linking of $S$ and $\partial Q$,  $g(Q)\cap S\neq
\varnothing$.

Set
$$
\Gamma(A)=\big\{m\in\C(\bar A,X);\ m=g \text{ on } \partial A\big\}
$$
The closure $\bar A$  of $A$ is compact and $X$ is a Banach 
space. Hence $\Gamma(A)$ endowed with the distance 
$$
\dist\nolimits_A(k_1,k_2)=\max_{x\in A}\| k_1(x),k_2(x)\|
$$
is a complete metric space.

Set
$$
\Psi(x)=\max\big\{0,\varepsilon^2-\varepsilon.\dist(x,S)\big\}
$$
and
$$
\I\colon \Gamma(A)\to \Rset,~ k\mapsto\I(k)=\max_{\bar A}\big\{f(k(t))+\Psi(f(t))\big\}
$$
$\bullet$ For each function $k\in \Gamma(A)$, we have
$$
k(\bar A)\cap S\neq \varnothing.
$$
Indeed, if we note by $\tilde k=k.\chi_{\overline A}+g.\chi_{A^c}$, 
it is obvious that $\tilde k\in E$ because $k=g$ on $\partial A$. Since $S$ 
and $ \partial Q$ link, there exists $\bar t\in Q$ such that $\tilde k 
(\bar t)\in S$. But by the definition of $A$,   
$\dist(g(t),S)\geq\varepsilon$ on $A^c\cap Q$, therefore $\bar t\in \bar A$.\\
$\bullet$ Using the relation ($*$), we have\footnote{Since $c=c_0=\inf_{S} f$
and $k(\bar t)\in S$.}
$$
\I(k)\geq f(k(\bar t))+\Psi(k(\bar t))\geq c+\varepsilon^2
$$
Since $k$ was taken arbitrary in $\Gamma(A)$,
\begin{equation}
\inf_{\Gamma(A)}\I\geq c+\varepsilon^2.
\end{equation}
But if we set $\tilde g=g|_{\bar A}$, we have
\begin{eqnarray}
\I(\tilde g)&=\max_{t\in \bar A}\{f(g(t))+\Psi(g(t))\},\\
            &\leq \max_{t\in Q}\{ f(g(t))+\Psi(g(t))\}.
\end{eqnarray}
Therefore
\begin{equation}
\I(\tilde g)\leq \left(c+\dfrac{\varepsilon^2}{4}\right)+\varepsilon^2
=c+\dfrac{5\varepsilon^2}{4}.
\end{equation}
Since $\I$ is bounded from below by $c+\varepsilon^2$ and lower semi-continuous
in the complete metric space $\Gamma(A)$ and $\tilde g\in \Gamma(A)$ verifies 
$$
\I(\tilde g)\leq \inf_{\Gamma(A)}\I+\dfrac{\varepsilon^2}{4}.
$$
Ekeland's variational principle implies then that there exists $\hat g\in \Gamma(A)$ such that 
\begin{equation}
\I(\hat g)\leq \I\tilde g),
\end{equation}
\begin{equation}
\| \hat g-\tilde g\|\leq \dfrac{\varepsilon}{2}
\end{equation}
and
\begin{equation}
\I(k)\geq\I(\hat g)-\dfrac{\varepsilon}{2}\|k-\hat g\|,\qquad \forall
k\in \Gamma(A)
\end{equation}
$\bullet$ Let $M$ be the subset of $\bar A$ where $f\circ \hat g$ realizes 
its maximum in $\bar A$.\\
\textsc{Claim:} \emph{There exists $t_0\in M$ such that:}
$$
\|f'(\hat g(t_0))\|\leq \dfrac{3}{2}\varepsilon.
$$
Indeed, for any $h\in \C(\bar A, X)$ such that $h|_{\partial A}\equiv 0$, 
we have
$$
f(\hat g(t_0)+\lambda h(t))=f(\hat g(t))+\lambda \langle f'(g(t)),h(t)\rangle 
+\circ(\lambda h(t))
$$
Therefore
$$
\max_{\bar A}f(\hat g(t_0)+\lambda h(t))\leq \max_{t\in \bar A}\{f(\hat g(t))+
\lambda \langle f'(g(t)),h(t)\rangle \}+\circ(\lambda\| h\|)
$$
where $\|h\| =\max_{\bar A}\|h(t))\|$.

Using this relation in (8), we obtain
$$
\I(g_\varepsilon)\geq \I(\hat g)+\dfrac{\varepsilon}{2}\|g_\varepsilon-\hat g\|
~~{\rm where}~~g_\varepsilon=\hat g+\lambda h.
$$
Therefore
$$
\max_{\bar A}f(\hat g(t))\leq \max_{t\in \bar A}\{f(\hat g(t))+
r\lambda \langle f'(g(t)),h(t)\rangle \}+{\varepsilon}\| h\|/{2}
$$
So denoting 
$$
\alpha(t)=f(\hat g(t)), \qquad \beta(t)=\langle f'(t),h(t)\rangle
$$
and
$$
N(\gamma)=\max_{\bar A }\gamma\quad\text{ with }\quad \gamma\in \C(\bar A,X).
$$
We have
$$
{N(\alpha+\lambda\beta)-N(\alpha)}/{\lambda}\geq -{\varepsilon} 
\|h\|/{2}, \quad \forall \lambda>0
$$
Therefore
\begin{equation}
\liminf_{\lambda\downarrow0}[N(\alpha+\lambda\beta)
-N(\alpha)]/{\lambda}\geq -{\varepsilon}\|h\|/{2}.
\end{equation}
 \begin{rk}
	 Notice that $\beta$ depends on $h$ and  that  $M\cap Q=\varnothing$.
\end{rk}

Set $N(\gamma)$ the subdifferential of $N$ at $\gamma$. We recall that 
$$
\partial N(\gamma)=\big\{\mu;\ \mu  {\text{ is a Radon measure of mass 1 with support
in }} M(\gamma)\big\} 
$$
Using Proposition 3 in~\cite{St}, we have:
 $$
\begin{array}{ll}
-\dfrac{\varepsilon}{2}\|h\|&\leq \liminf\limits_{\lambda\downarrow0}[N(\alpha+\lambda
\beta)-N(\alpha)] \\ 
&\leq \max\big\{\langle \beta,\mu\rangle / \mu \in \partial N(\alpha)\big\}\\
&\ \ =\max\left\{\dint_{\bar A} \langle f'(\hat g), h\rangle\,d\mu;\ \mu \in\partial
N(\alpha)\right\} 
\end{array}   
$$
Now, apply \cite[Theorem~6.2.7]{AE}  to the functional:
$$
\mathcal{G}\colon\mathcal{M}(\bar A,\Rset)\times\C(A,X)\to\Rset
$$
$$
\mathcal{G}(\mu,k)=\langle \mu,\langle f'(\hat g(.)),k(t)\rangle\rangle=\dint_{\bar
A}\langle f'(\hat g),h\rangle\,d\mu
$$
Let $\mathcal{M}(K,\Rset)$ be the Banach space of Radon measures defined on the
$\sigma$-Algebra of all Borel sets of $\bar A$ endowed with the
$w^*$-topology. Then $\mathcal{G}$ is continuous, linear in each variable 
separately. And the sets $\partial N(\alpha)$ and $\big\{ k\in 
\C(K,X);\ \|k\|\leq  1\big\}$  are convex, the former one 
being $w^*$-compact. We have    
$$
\begin{array}{ll}
-\dfrac{\varepsilon}{2}\|h\|&\leq \inf_{h}\max_{\mu}\left\{\dint_{\bar A} \langle f'(\hat
g), h\rangle\,d\mu;\ \mu \in\partial N(\alpha), \|h\|\leq1,h|_{\partial 
\bar A}=0\right\} \\
&\ \ =\max_{\mu}\inf_{h}\left\{\dint_{\bar A} \langle f'(\hat
g), h\rangle\,d\mu;\ \mu \in\partial N(\alpha), \|h\|\leq1,h|_{\partial \bar A}=0\right\} \\
&\ \ =\max_{\mu}\left\{\dint_{\bar A} \| f'(\hat g) \|\;\ \mu \in\partial N(\alpha)\right\} \\
&\ \ =-\min\{ \| f'(\hat g(t) )\|;\ t\in M(\underbrace{f'(\hat g)}_{\alpha})\} \\
\end{array}
$$
Therefore, there exists $t_0\in M$ such that $\|f'(\hat g(t_0))\|\leq  \dfrac{3}{2}
\varepsilon $ and the Claim is proved.\\
\textbf{For (i):} By (4), (5), (6) and the Claim, we have  
$$
c+\varepsilon^2\leq\inf_{\Gamma(A)}\I \leq f(\hat g(t_0)+\Psi(\hat g(t_0))=\I(\hat
g)\leq\I(\tilde g)\leq c+\dfrac{5}{4}\varepsilon^2.
$$
Since $0\leq\Psi\leq\varepsilon^2$, we have $c\leq f(x_\varepsilon)\leq c+
\dfrac{5}{4}\varepsilon^2$.\\
\textbf{ For (ii):} It suffices to see that since $t_0\in \bar A= \big\{ 
x\in Q;\ \dist(g(x),S) \leq \varepsilon\big\}$, we have $\dist(\tilde 
g(t_0),S)=\dist(g(t_0),S)\leq  \varepsilon$. Then by (7), we  
have 
$$
\dist(x_\varepsilon,S)=\dist(\hat g(t_0),S)\leq  \dfrac{3}{2}
$$
\cqfd       

\begin{rk}
	Assuming that $f$ verifies (P.S.)$_c$, the inf max value $c$ is critical 
	in the limiting case and there exits a critical point of level $c$ in 
	$S$. This confirms the result obtained by the deformation lemma.
\end{rk}
\begin{rk}
	In~\cite{Ma}, Mawhin has pointed out that in such situations we need a weaker
	condition than (P.S.)$_c$ to conclude. Indeed, since we are sure that there exists a
	sequence $(u_k)$ such that:
	$$
	f(u_k)\to c  \quad\text{and }\quad f'(u_k)\to 0
	$$
	It suffices then to assume that $f$ is such that the existence of such a sequence
	implies that $c$ is a critical value.

	Unfortunately, in applications, $c$ isn't known explicitly. This constrains us to
	verify (P.S.)$_c$ for all possible values $c$.
\end{rk}

As corollaries of the abstract results, we obtain the known theorems:
\begin{cor}[Rabinowitz~\cite{Ra}]
Let $\Phi\in \C^1(X,\Rset)$ a functional that satisfies \text{(P.S.)}. Suppose that 
$$
\inf_{u\in B(0,r)}\Phi(u)=\max\{\Phi(0),\Phi(e)\}=c=\inf_{\gamma\in\Gamma}
\max_{t\in[0,1]}\Phi(\gamma(t))
$$
where $0<r<\|e\|$ and $\Gamma=\big\{\gamma\in\C([0,1],X);\  
\gamma(0)=0 \text{ and }\gamma(1)=e\big\}$\\
Then $\Phi$ has a critical point of level $c$ on the sphere $S(0,r)$.
\end{cor}
\begin{cor}[Pucci-Serrin~\cite{P-S}]
Let $\Phi\in \C^1(X,\Rset)$ satisfying {\rm (P.S.)}, if $\Phi$ has a pair of local minima
(or local maxima) then $\Phi$ possesses a third critical point.
\end{cor}

 \section*{APPENDIX}
 In the proof of the linking in Examples 1 and 2, we use the homotopy property 
 of the degree theory respectively with: 
 $$
  F_t \colon Q \to V_2, \quad\text{where}\quad F_t (u)=tP_2 \gamma (u) + 
 (1-t)u \qquad \qquad t \in \left[0,1\right]
 $$  
 and to
 $$
  F_t \colon Q \to \overline {\Rset \times V_2}, \quad {\hbox {\rm such that}}
 \quad \text{for all } u \in V_2 \cap B_{R_2} \quad\text{and for all }s 
 \in \left[0,R_1\right]
 $$  
 $$ 
 F_t (u)=\big( {\left|| tP_1 \gamma (u+se) \right||} - \rho + (1-t)s, 
 tP_2 \gamma (u+se)  + (1-t)u \big) 
 $$
 It is obvious there that the finite dimension of $V_2$ can be dropped 
 if we require $\gamma$ to be 
 compact and the degree will still have a sense as  
 Schauder's degree of a compact perturbation of the identity.
 \vskip 2mm
 \noindent \textbf{Proof of the linking of $Q_R$ and $S_{\rho}$ in 
 Theorem~\protect\ref{silva2}. } 

 We want to prove it without assuming $\dim V_2<\infty$. When assuming
$\dim V_2\leq \infty$ the proof becomes easier. We will adapt the proof of Silva to our
situation. The compactness of the elements  of
$\Gamma_4$, $\Gamma_5$, $\Gamma_2$, $\Gamma_2 \cap \Gamma_3$ compensating
the local compactness of the finite dimensional space and allowing the use of
 degree theory as explained above.\\
The space $V=V_1 \oplus \Rset e \oplus V_2$, while the reals $\rho$ 
and $R$ are such that $0<\rho<R$ $D_R = \partial B_R(0) \cap
 (V_1 \oplus \Rset e )$ and $S_{\rho} = \overline{ V_2 \setminus B_{\rho} (0) )}
 \cup \left( \partial B_{\rho}(0) \cap (\Rset^+ e \oplus V_2\right)$, \\
 If we take $\Gamma =\Gamma_4$ or $\Gamma =\Gamma_5$ or $\Gamma =\Gamma_2 \cap 
 \Gamma_3$ \ then $S_{\rho}$\ and\ $D_R$ link with respect to~ 
 $\Gamma$. Indeed: \\
 
 \noindent $\bullet $ First, remark that $S_{\rho} \cap \partial D_R = \varnothing$  
 because if $x \in S_{\rho} \cap \partial D_R$, then we would have 
 \vskip 1mm
 \begin{center}
 $\left\{~ \begin{tabular}{ll}
         $|| x || = R$ &$ \Rightarrow || x || = R$,\\
       $  \left.\begin{tabular}{l}
         $\! \! \! x \in V_1 \oplus \Rset e $\\
         $\! \! \! x \in V_2$
         \end{tabular}
         \right\}$
                 &$ \Rightarrow x  = 0. $
        \end{tabular}
        \right.$
 \end{center}
 And this yields a contradiction.\\
 
 \noindent $\bullet $ We claim  that for all $\gamma \in \Gamma$, it 
 holds that $\gamma(D_R)\cap S_{\rho} \neq \varnothing$. Indeed, let
 $$ 
 \chi_{\beta} =\left\{~  
 \begin{tabular}{lll}
 1&if&$x \geq  \rho / \beta$ \\
 $\beta x / \rho$ &if &$ x \in [0,\rho / \beta ]$\\
 0& if & $x \leq  0 $
 \end{tabular}
 \right. $$ 
 Notice that $\chi_{\beta} \in {\C}(\Rset,\Rset)$, let $P_1$ be the 
 projection $P_1 \colon V \to V_1$ and
 $P_e \colon V \to \Rset e$ the projection onto $\Rset e$. Set 
 $$ 
 \begin{tabular}{ll}
 $G_{\beta,s}(re+v) =$&$ \Biggl( s(re+v)  +(1-s) \Biggl[ P_1 
 \biggl( \gamma (re+v) \biggr) +        $ \\
  ~&$~~~~~~ \biggl[ \chi_{\beta} \biggl(P_e  \bigr(\gamma (re +v) \bigr) \biggr) 
 \Big|| (I_d -P_1) \bigl( \gamma (re+v) \bigr) \Big|| \biggr] e 
 \Biggr] \Biggr) $
  \end{tabular}
  \eqno(1) 
  $$
 where $re+v$ denotes a generic element of $V_1 \oplus \Rset e$  and $s \in [0,1]$, 
 then $G_{\beta,s} \in {\C}(V_1 \oplus \Rset e , V_1 \oplus \Rset e)$.\\ 
 If $re + v \in \partial D_R $, since $\gamma |_{\partial D_R}=I_d$,  
 we have 
 $$ G_{\beta,s}(re+v) = \bigg(sr + v + (1-s) \chi_{\beta} (r) | r| \bigg) 
 \eqno(2)$$
 Therefore, if $\quad\beta >1$ we have that $\rho e \not\in G_{\beta,s}(\partial D_R)$
 for any $s \in [0,1]$ because if $ G_{\beta,s}(re+v) =\rho e$, then $r= R$
 and $v=0$, and this implies that $ G_{\beta,s}(re+v) =R e$, and this 
 contradicts the fact that $\rho < R$.\\
 So, the topological degree of the compact perturbation of the identity $G_{\beta,s}$
 on $D_R$ at $\rho e$ is well defined. By the homotopy invariance of 
 the degree, we have
 $$ d(G_{\beta,s},D_R,\rho e)=d(I_d,D_R, \rho e)=1,$$
 Therefore, for $s=0$ we conclude that  
 $$
 \left\{~  \begin{tabular}{l}
           $ P_1 \gamma (u_{\beta}) = 0 $, \\
           $ \chi_{\beta} (P_e \gamma(u_{\beta})) \cdot || (I_d -P_1)\gamma 
           (u_{\beta})|| = \rho >0 $.
           \end{tabular}
 \right. \eqno(3) 
 $$
 By the definition of $\chi_{\beta}$ and the second equation of (3), we 
 have that $P_e \gamma (u_{\beta}) > 0$. Then, if we suppose that 
 $\gamma(D_R) \cap S_{\rho} = \varnothing$ we would have $ P_e 
 (\gamma(u_{\beta})) < \rho / \beta $ where   
 $\beta >1$. If that wasn't the case, $\chi_{\beta}\big(P_e \gamma
 (u_{\beta})\big)=1$ and hence $|| (I_d -P_1) \gamma (u_{\beta})|| =\rho$. But
 since  $P_1 \big( \gamma (u_{\beta})\big)=0$, we would obtain $|| \gamma
 (u_{\beta})   || = \rho$ and $ \gamma  (u_{\beta}) \in V_2 \oplus \Rset e, \quad$ i.e.
 $\gamma (u_{\beta}) \in S_{\rho}$? A contradiction with what we supposed before, so
 $P_e \gamma (u_{\beta}) < \rho /  \beta$.\\
 Let ${(\beta_m)}_m \subset \Rset $ a sequence such that $\beta_m \to \infty $ as 
 $m \to \infty $ and  ${(u_{\beta_m})}_m \subset D_R $ the sequence given by 
 (3) and satisfying 
 $$
 0< P_e (\gamma (u_{(\beta_n)}) < \rho / \beta_n
 \eqno(4)
 $$
 Since $D_R$ is bounded and ${(u_{(\beta_n)})}_n \subset D_R$, the 
 sequence  ${(\gamma 
 (u_{(\beta_n)}))}_n$ is relatively compact. So, for a subsequence still denoted 
 by $ {(u_{(\beta_n)})}_n$ we have $\gamma (u_{(\beta_n)}) \to u \in V_2 
 \oplus \Rset e $ because $\gamma \in \Gamma$ is compact.\\
 By (4) $P_e u = \lim_n P_e \gamma(u_{(\beta_n)}) =0$ and since $P_1 
 \gamma(u_{(\beta_n)})=0$   for any $n \in \Nset$ we have    
 $$ 
 P_e u = P_1 u = 0.
 $$
 But $\chi_{\beta} \leq 1$, so by the second equation of (3) we have  
 $$
 || u || = || (I_d - P_1) u || = \lim_n || (I_d - P_1) 
 \gamma(u_{(\beta_n)}) \| \geq \rho  
 $$ 
 and hence $u \in \overline{\big(V_2 \setminus B_{\rho} (0) \big)} \subset S_{\rho} $
 i.e. $ u \in \overline{ \gamma (D_R )} \cap S_{\rho}$. While by the definition of 
  $\Gamma$ taken here  we have $\overline{ \gamma (D_R )}= 
  \gamma (D_R )$. A contradiction with $\gamma (D_R ) \cap S_{\rho} = \varnothing$.
  \cqfd

 \begin{rk}
	 We saw that the finite dimension of $V_2$ isn't required to prove the linking in
	any of the cases presented here when we use the smaller sets of ``compact''
	functions, so why can't we use them to remove this additional 
	assumption? The answer is that
	these sets are so small that they are empty. At least this happens with Hilbert
	spaces, indeed,
	let $H=V_1\oplus V_2$ with both $\dim V_1$ and $\dim V_2$ infinite, and suppose by
	contradiction that the set  
	$$
	\Gamma=\big\{\gamma\in \C(H,H)\text{ compact};\ \gamma|_{S(0,R)\cap V_2}=I_d\big\}
	$$
	is not empty. Without loss of generality, we can suppose $R=1$ so that $S(0,R)$ is the
	unit sphere. Let us take an orthonormal basis $(u_n)_{n}$ of $V_2$, we know that it
	converges weakly to 0. Take $\gamma\in\Gamma$, since the unit sphere is
	bounded and $\gamma$ compact, $\gamma(S(0,R))$ is compact and the sequence
	$\gamma(u_n)=u_n$ admits a subsequence converging to an element of $S(0,R)$,
	 a contradiction with the fact that $u_n \rightharpoonup 0$.
\end{rk}
 
 \noindent \textbf{Proof of Lemma~\ref{lemme1}\/} \\
 The details will be omitted since the proof is not difficult and classical 
 but quite long. 
 It suffices to change $f'$ by a pseudogradient field in the proof of~\cite{Y}.
 We sketch it here in few lines for the convenience of the reader.\\
 The map $\eta$ is obtained as a solution of a differential equation. First we 
 easily verify that \\
 ${\rm (1)} \overline{N_{\delta}(E)} \cap \mathcal{K}_c = \varnothing$, and\\
 ${\rm (2)} \left|| f'(x) \right|| \geq b$,  for all $x \in 
 A=N_{\delta}(E) \cap f_{c-{\bar \varepsilon}}^{c+{\bar \varepsilon}}$.\\
 Set $A_1 = V \setminus  \left( \overline {N_{\delta / 2}(E)} \cap f_{c-{\bar 
 \varepsilon / 2}}^{c+{\bar \varepsilon / 2}} \right)$, and   $A_2=
 N_{\delta / 3}(E) \cap f_{c-{\bar \varepsilon / 3}}^{c+{\bar 
 \varepsilon / 3}}$  \\
 and define 
 $$
 \begin{tabular}{rl} 
 \\ $h(x) =$ &${ {|| x- A_1 ||} / { \left( || x- A_1
||  + || x- A_2 || \right)}}, $ \cr
 \\ 
 $ \rho (s) =$ & $ \left\{~ \begin{tabular}{rc}
                               $1  \text{if }$ & $0 \leq s \leq 1$ \cr
                               $1 / s \text{ if }$ & $ s \geq 1$\cr
                           \end{tabular} \right.$ \\
 \\ $ X_1 (x) =$ & $ \left\{~        
 \begin{tabular}{ll}
$-\delta /3 ~h(x) \rho \left( || W(x) || \right) W(x) $ &
$\text{ in }A$ \cr
$0 $ & $ \text{ in }V \setminus A $\cr
\end{tabular}   \right.$
 \end{tabular}
 $$
 \\where $W$ is a pseudogradient field defined on $\tilde V$. And let 
 $\zeta $ be the solution of 
 $$ \left\{~ \begin{tabular}{lll}
               $ {d\zeta} /{ dt}$& =&   $ X_1 ( \zeta ) $\cr
                               $\zeta (0,x)$& =&  $ x \in V$\cr
              \end{tabular} \right.
$$
 Then, if $E_1 =\zeta \left( [0,1]  \times E \right)$, we prove easily using the graph norm
 that $E_1$ is closed and with some obvious properties of $\zeta$ that $E_1 
 \cap D = \varnothing$. \\
 Let $D_1 =\left\{ x \in V;\ d(x) \leq 1\right\}$, where 
 $ d(x) = { {\left|| x-D \right|| /  \left( \left|| x-D \right||
 + \left|| x - E_1 \right|| \right)}} $. 
 It is clear then that $D_1 \cap E_1 = \varnothing$.\\
 Define now $X(x) =  q(x) X_1(x)$ with $q(x)= { {\left|| 
 x-D_1 \right|| /  \left( \left|| x-D_1 \right|| 
 + \left|| x-E_1 \right|| \right)} }$ 
 then the solution $\eta$ of
 $$\left\{~ \begin{tabular}{lll}
            $ { d \eta} /{ dt} $& =&  $X ( \eta ) $ \\
            $ \eta (0,x) $&= & $ x \in V $ \\
           \end{tabular}
  \right.
  $$
 satisfies the desired properties (the proof there is standard except 
 for some  obvious changes in (iv) ).
\newline 
\textbf{Note.} Definition 1 has its origin in~\cite{St}. After a first version of
this paper was written, it was brought to our attention that the same definition appears
also in~\cite{Wi}. In Willem's paper an
interesting result near to our Theorem~\ref{deuxieme} is given~\cite[Lemme 2.18]{Wi}:
\begin{thm}[Th\'eor\`eme de localisation de Willem]
\label{Willem}
Let $X$ be a Banach space, $\Phi\in\mathcal{C}^1(X,\Rset)$, $F,Q\subset X$. Write
$$
   c=\inf_{\gamma\in \Gamma}\sup_{u\in Q}\Phi(\gamma(u))
$$
where
$$
\Gamma=\big\{\gamma\in\mathcal{C}(Q,X);\ \gamma(u)=u \text{ on }\partial Q\big\}
$$
If\\
(D1) $F$ and $\partial Q$ link,\\
(D2) dist$(F,\partial Q) >0$,\\
(D3) $-\infty<c=\inf_{F}\Phi$,\\
Then for each $\varepsilon >0$, $\delta\in]0,\dist(F,\partial Q)/2[$ and $\gamma\in
\Gamma$ such that
$$
\sup_{Q}\Phi\circ\gamma<c+\varepsilon
$$
There exists $u\in X$ such that\\
a) $c-2\varepsilon\leq\Phi(u)\leq c+2\varepsilon$\\
b) dist$(u,F\cap\gamma(Q)_\delta)\leq2\delta$\\
c) $\|\Phi'(u)\|<4\varepsilon/\delta$.
\end{thm}
We report some remarks
\begin{itemize}
\item In Willem's result, $Q$ isn't supposed compact so the assumption $-\infty <c$ which
appears in (D3) is necessary.
\item When $Q$ is compact, our result contains Theorem~\ref{Willem}. 
Our result contains also classical point theorems cited above while this 
isn't the case for~\ref{Willem}.  
\item The proof of our Theorem~\ref{deuxieme} uses ``Ekeland's Variational 
principle'' while Willem 
uses a \emph{General Deformation Lemma} to prove his result. 
\end{itemize}
{\small 
{\sl Acknowledgment}---We are also
 grateful to  Professor M.~Willem  for sending us his reference~\cite{Wi}.
 }

\end{document}